\documentclass[11pt]{article}
\usepackage{amssymb,amsmath,graphicx,theorem}

\def\GG{{\cal G}}

\def\KK{{\cal K}}

\def\PP{{\cal P}}

\def\R{{\mathbf R}}
\def\Z{{\mathbf Z}}

\def\T{{^{\sf T}}}

\def\hom{\hbox{\rm hom}}

\def\rk{\hbox{\rm rk}}

\newtheorem{theorem}{Theorem}[section]
\newtheorem{prop}[theorem]{Proposition}

\newtheorem{claim}{Claim}[section]

\theorembodyfont{\rmfamily}

\newtheorem{example}[theorem]{Example}

\newenvironment{proof}{\noindent{\bf Proof. }}{\hfill$\square$\medskip}

\begin{document}

\title{Reflection positivity, rank connectivity, and homomorphism of graphs
}
\author{{\sc Michael Freedman}, {\sc L\'aszl\'o Lov\'asz}\\
Microsoft Research \\
One Microsoft Way\\
Redmond, WA 98052\\
and
\\
{\sc Alexander Schrijver}\\
CWI Amsterdam\\
The Netherlands}

\date{March 2004\\[1cm] MSR-TR-2004-41}

\maketitle

\begin{abstract}
It is shown that a graph parameter can be realized as the number of
homomorphisms into a fixed (weighted) graph if and only if it
satisfies two linear algebraic conditions: reflection positivity and
exponential rank-connectivity. In terms of statistical physics, this
can be viewed as a characterization of partition functions of vertex
models.
\end{abstract}

\section{Introduction}

For two finite graphs $G$ and $H$, let $\hom(G,H)$ denote the number
of homomorphisms (adjacency-preserving mappings) from $G$ to $H$.
Many interesting graph parameters can be expressed in terms of the
number of homomorphisms into a fixed graph: for example, the number
of colorings with a number of colors is the number of homomorphisms
into a complete graph. Further examples for the occurrence of these
numbers in graph theory will be discussed in Section \ref{EXAMPLES}.

Another source of important examples is statistical physics, where
partition functions of various models can be expressed as graph
homomorphism functions. For example, let $G$ be a $n\times n$ grid,
and suppose that every node of $G$ (every ``site'') can be in one of
two states, ``UP'' or ``DOWN''. The properties of the system are such
that no two adjacent sites can be ``UP''. A ``configuration'' is a
valid assignment of states to each node. The number of configurations
is the number of independent sets of nodes in $G$, which in turn can
be expressed as the number of homomorphisms of $G$ into the graph $H$
consisting of two nodes, "UP" and "DOWN", connected by an edge, and
with an additional loop at "DOWN". To capture more interesting
physical models, so called ``vertex models'', one needs to extend the
notion of graph homomorphism to the case when the nodes and edges of
$H$ have weights (see Section \ref{WEIGHT}).

Which graph parameters can be represented as homomorphism functions
into weighted graphs? This question is motivated, among others, by
the problem of physical realizability of certain graph parameters.
Two necessary conditions are easy to prove:

\smallskip

(a) The interaction between two parts of a graph separated by $k$
nodes is bounded by a simple exponential function of $k$ (this will
be formalized as {\it rank-connectivity} in Section \ref{RANKCON}).

\smallskip

(b) Another necessary condition, which comes from statistical
mechanics, as well as from extremal graph theory, is {\it reflection
positivity}. Informally, this means that if a system has a 2-fold
symmetry, then its partition function is positive. We'll formulate a
version of this in section \ref{RANKCON}.

\smallskip

The main result of this paper is to prove that these two necessary
conditions are also sufficient. The proof makes use of a simple kind
of (commutative, finite dimensional) $C^*$-algebras.

\section{Homomorphisms, rank-connectivity and reflection positivity}

\subsection{Weighted graph homomorphisms}\label{WEIGHT}

A {\it weighted graph} $H$ is a graph with a positive real weight
$\alpha_H(i)$ associated with each node $i$ and a real weight
$\beta_H(i,j)$ associated with each edge $ij$.

Let $G$ be an unweighted graph (possibly with multiple edges, but no
loops) and $H$, a weighted graph. To every homomorphism
$\phi:~V(G)\to V(H)$, we assign the weight
\begin{equation}\label{HOMW}
w(\phi)=\prod_{u\in V(G)} \alpha_H(\phi(u)) \prod_{uv\in E(G)}
\beta_H(\phi(u),\phi(v)),
\end{equation}
and define
\begin{equation}\label{HOMDEF}
\hom(G,H)=\sum_{\phi:~V(G)\to V(H)\atop \text{homomorphism}} w(\phi).
\end{equation}
If all the node-weights and edge-weights in $H$ are $1$, then this is
the number of homomorphisms from $G$ into $H$ (with no weights).

For the purpose of this paper, it will be convenient to assume that
$H$ is a complete graph with a loop at all nodes (the missing edges
can be added with weight $0$). Then the weighted graph $H$ is
completely described by a positive integer $d=|V(H)|$, the positive
real vector $a=(\alpha_1,\dots,\alpha_d)\in\R^d$ and the real
symmetric matrix $B=(\beta_{ij})\in\R^{d\times d}$. The graph
parameter $\hom(.,H)$ will be denoted by $f_H$ or $f_{B,a}$.

\subsection{Connection matrices of a graph parameter}\label{RANKCON}

A {\it graph parameter} is a function on finite graphs (invariant
under graph isomorphism). We allow multiple edges in our graphs, but
no loops. A graph parameter $f$ is called {\it multiplicative}, if
for the disjoint union $G_1\cup G_2$ of two graphs $G_1$ and $G_2$,
we have $f(G_1\cup G_2)=f(G_1)f(G_2)$.

A {\it $k$-labeled graph} ($k\ge 0$) is a finite graph in which $k$
nodes are labeled by $1,2,\dots k$ (the graph can have any number of
unlabeled nodes). Two $k$-labeled graphs are {\it isomorphic}, if
there is a label-preserving isomorphism between them. We denote by
$K_k$ the $k$-labeled complete graph on $k$ nodes, and by $O_k$, the
$k$-labeled graph on $k$-nodes with no edges. In particular,
$K_0=O_0$ is the graph with no nodes and no edges.

Let $G_1$ and $G_2$ be two $k$-labeled graphs. Their {\it product}
$G_1G_2$ is defined as follows: we take their disjoint union, and
then identify nodes with the same label. Hence for two $0$-labeled
graphs, $G_1G_2=G_1\cup G_2$ (disjoint union).

Now we come to the construction that is crucial for our treatment.
Let $f$ be any graph parameter. For every integer $k\ge 0$, we define
the following (infinite) matrix $M(f,k)$. The rows and columns are
indexed by isomorphism types of $k$-labeled graphs. The entry in the
intersection of the row corresponding to $G_1$ and the column
corresponding to $G_2$ is $f(G_1G_2)$. We call the matrices $M(f,k)$
the {\it connection matrices} of the graph parameter $f$.

We will be concerned with two properties of connection matrices,
namely their rank and positive semidefiniteness. The rank
$r_f(k)=\rk(M(f,k))$, as a function of $k$, will be called the {\it
rank connectivity function} of the parameter $f$. This may be
infinite, but for many interesting parameters it is finite, and its
growth rate will be important for us. We say that a set of
$k$-labeled graphs {\it spans} if the corresponding rows of $M(f,k)$
have rank $r_f(k)$.

We say that a graph parameter $f$ is {\it reflection positive}, if
$M(f,k)$ is positive semidefinite for every $k$. This property is
closely related to the reflection positivity property of certain
statistical physical models. Indeed, it implies that for any
$k$-labeled graph $G$, $f(GG)\ge 0$ (since $f(GG)$ is a diagonal
entry of a positive semidefinite matrix). Here the second copy of $G$
can be thought of as a ``reflection'' (in the set of labeled nodes)
of the first.

The condition that $f(GG)\ge 0$ is weaker than the condition that
$M(f,k)$ is positive semidefinite, but we can strengthen it as
follows. Let $x=(x_G)$ be a complex vector indexed by $k$-labeled
graphs with $\sum_G |x_G|^2=1$. The formal linear combination
$X=\sum_G x_G G$ can be thought of as a ``quantum $k$-labeled
graph'', and then $X\overline{X}=\sum_{G_1,G_2}
x_{G_1}\overline{x}_{G_2} (G_1G_2)$ is a quantum graph that can be
obtained by gluing together $X$ and its reflection. (Note that we
took complex conjugation of the coefficient to build the ``mirror
copy'' of $X$.) Now if we extend $f$ linearly over quantum graphs,
then
\[
f(XX) = x\T M(f,k) \overline{x},
\]
showing that the non-negativity of $f$ over symmetric quantum graphs
is equivalent to reflection positivity as we defined it.

\subsection{Simple properties of connection matrices}

We give a couple of simple facts about connection matrices of a
general graph parameter.

\begin{prop}\label{MULT}
Let $f$ be a graph parameter that is not identically $0$. Then $f$ is
multiplicative if and only if $M(f,0)$ is positive semidefinite,
$f(K_0)=1$, and $r_f(0)= 1$.
\end{prop}

\begin{proof}
If $f$ is multiplicative, then $f(K_0)^2=f(K_0)$, showing that
$f(K_0)\in\{0,1\}$. If $f(K_0)=0$, then the relation
$f(G)=f(GK_0)=f(G)f(K_0)$ implies that $f(G)=0$ for every $G$, which
was excluded. So $f(K_0)=1$. Furthermore, $f(G_1G_2)=f(G_1)f(G_2)$
for any two $0$-labeled graphs $G_1$ and $G_2$, which implies that
$M(f,0)$ has rank 1 and is positive semidefinite.

Conversely, since $M(f,0)$ is symmetric, the assumption that
$r_f(0)=1$ implies that there is a graph parameter $\phi$ and a
constant $c$ such that $f(G_1G_2)=c\phi(G_1)\phi(G_2)$. Since
$M(f,0)$ is positive semidefinite, we have $c > 0$, and so we can
normalize $\phi$ to make $c=1$. Then $f(K_0)=f(K_0K_0)=\phi(K_0)^2$,
whence $\phi(K_0)\in\{-1,1\}$. We can replace $\phi$ by $-\phi$, so
we may assume that $\phi(K_0)=1$. Then
$f(G)=f(GK_0)=\phi(G)\phi(K_0)=\phi(G)$ for every $G$, which shows
that $f$ is multiplicative.
\end{proof}

\begin{prop}\label{SUPERMULT}
Let $f$ be a multiplicative graph parameter and $k,l\ge 0$. Then
\[
r_f(k+l)\ge r_f(k)\cdot r_f(l).
\]
\end{prop}

We'll see (see Claim \ref{RESOLVEDE(G)} below) that in the case when
$f$ is reflection positive, the sequence $r_f(k)$ logconvex. We don't
know if this property holds in general.

\medskip

\begin{proof}
Let us call a $(k+l)$-labeled graph {\it separated}, if every
component of it contains either only nodes with label at most $k$, or
only nodes with label larger than $k$. Consider the submatrix of
$M(f,k+l)$ formed by the separated rows and columns. By
multiplicativity, this submatrix is the Kronecker (tensor) product of
$M(f,k)$ and $M(f,l)$, so its rank is $r_f(k)\cdot r_f(l)$.
\end{proof}

\subsection{Connection matrices of homomorphisms}

Fix a weighted graph $H=(a,B)$. For any $k$-labeled graph $G$ and
mapping $\phi:~[1,k]\to V(H)$, let
\begin{equation}\label{PHI-HOM}
\hom_\phi(G,H)=\sum_{\psi:~V(G)\to V(H)\atop \psi
\text{~extends~}\phi} w(\psi).
\end{equation}
So $\hom(G,H)=\sum_\phi \hom_\phi(G,H)$.

\begin{theorem}\label{EASY}
The graph parameter $f_H$ is reflection positive and $r_{f_H}(k)\le
|V(H)|^k$.
\end{theorem}

\begin{proof}
For any two $k$-labeled graph $G_1$ and $G_2$ and $\phi:[1,k]\to
V(H)$,
\begin{equation}\label{PHI-PROD}
\hom_\phi(G_1G_2,H)=\frac{1}{\prod_{i=1}^k \alpha_H(\phi(i))}
\hom_\phi(G_1,H)\hom_\phi(G_2,H).
\end{equation}
The decomposition (\ref{PHI-HOM}) writes the matrix $M(f,k)$ as the
sum of $|V(H)|^k$ matrices, one for each mapping $\phi:~[1,k]\to
V(H)$; (\ref{PHI-PROD}) shows that these matrices are positive
semidefinite and have rank 1.
\end{proof}

The main result of this paper is a converse to Theorem \ref{EASY}:

\begin{theorem}\label{MAIN}
Let $f$ be a reflection positive graph parameter for which there
exists a positive integer $q$ such that $r_f(k) \le q^k$ for every
$k\ge 0$. Then there exists a weighted graph $H$ such that $f=f_H$.
\end{theorem}

\section{Examples}\label{EXAMPLES}

We start with an example showing that exponential growth of rank
connectivity is not sufficient in itself to guarantee that a graph
parameter is a homomorphism function.

\begin{example}[Matchings]\label{MATCHING}
Let $\Phi(G)$ denote the number of perfect matchings in the graph
$G$. It is trivial that $\Phi(G)$ is multiplicative. We claim that
its node-rank-connectivity is exponentially bounded:
\[
r_\Phi(k)=2^k.
\]

Let $G$ be a $k$-labeled graph, let $X\subseteq [1,k]$, and let
$\Phi(G,X)$ denote the number of matchings in $G$ that match all the
unlabeled nodes and the nodes with label in $X$, but not any of the
other labeled nodes. Then we have for any two $k$-labeled graphs
$G_1,G_2$
\[
\Phi(G_1G_2) =\sum_{X_1\cap X_2=\emptyset,\,X_1\cup X_2=[1,k]}
\Phi(G_1,X_1)\Phi(G_2,X_2).
\]
This can be read as follows: The matrix $M(\Phi,k)$ can be written as
a product $N\T W N$, where $N$ has infinitely many rows indexed by
$k$-labeled graphs, but only $2^k$ columns, indexed by subsets of
$[1,k]$,
\[
N_{G,X}=\Phi(G,X),
\]
and $W$ is a symmetric $2^k\times 2^k$ matrix, where
\[
W_{X_1,X_2}=
  \begin{cases}
    1 & \text{if $X_1=[1,k]\setminus  X_2$}, \\
    0 & \text{otherwise}.
  \end{cases}
\]
Hence the rank of $M(\Phi,k)$ is at most $2^k$ (it is not hard to see
that in fact equality holds).

On the other hand, let us consider $K_1$ and $K_2$ as 1-labeled
graphs. Then the submatrix of $M(\Phi,1)$ indexed by $K_1$ and $K_2$
is
\[
\begin{pmatrix}0&1\\1&0\end{pmatrix}
\]
which is not positive semidefinite. Thus $\Phi(G)$ cannot be
represented as a homomorphism function.
\end{example}

\begin{example}[Flows]\label{FLOWS}
The following example shows that there are important graph parameters
that are not defined as homomorphism functions, but that can also be
represented as homomorphism functions in a nontrivial way. In fact,
it is easier to check that the conditions in Theorem \ref{MAIN} hold
than to show that these graph parameters are homomorphism functions,
and they first came up as counterexample candidates.

Let us start with a simple special case. Let $f(G)=1$ if $G$ is
eulerian (i.e., all nodes have even degree), and $f(G)=0$ otherwise.
To represent this function as a homomorphism function, let
\[
a=\begin{pmatrix}1/2\\1/2\end{pmatrix},\qquad
B=\begin{pmatrix}1&-1\\-1&1\end{pmatrix}
\]
It is not hard to see that for the weighted graph $H=(a,B)$ we have
$\hom(G,H)=f(G)$.

It follows that for this function, reflection positivity holds, and
the rank-connectivity is at most $2^k$.

This example can be generalized quite a bit. Let $\Gamma$ be a finite
abelian group and let $S\subseteq\Gamma$ be such that $S$ is closed
under inversion. For any graph $G$, fix an orientation of the edges.
An {\em $S$-flow} is an assignment of an element of $S$ to each edge
such that for each node $v$, the product of elements assigned to
edges entering $v$ is the same as the product of elements assigned to
the edges leaving $v$. Let $f(G)$ be the number of $S$-flows. This
number is independent of the orientation.

The choice $\Gamma=\Z_2$ and $S=\Z_2\setminus \{0\}$ gives the
special case above (incidence function of eulerian graphs). If
$\Gamma=S=\Z_2$, then $f(G)$ is the number of eulerian subgraphs of
$G$. Perhaps the most interesting special case is when $|\Gamma|=t$
and $S=\Gamma\setminus\{0\}$, which gives the number of nowhere zero
$t$-flows.

Surprisingly, this parameter can be described as a homomorphism
function. Let $\Gamma^*$ be the character group of $\Gamma$. Let $H$
be the complete directed graph (with all loops) on $\Gamma^*$. Let
$\alpha_{\chi}:=1/|\Gamma|$ for each $\chi\in\Gamma^*$, and let
\[
\beta_{\chi,\chi'}:=\sum_{s\in S}\chi^{-1}(s)\chi'(s),
\]
for $\chi,\chi'\in\Gamma^*$.

We show that
\begin{equation}\label{FLOWHOM}
f=\hom(.,H).
\end{equation}
Let $n=|V(G)|$ and $m=|\Gamma|$. For any coloring $\phi:~V(G)\to S$
and node $v\in V(G)$, let
\[
\partial_\phi(v)=\sum_{u\in V(G)\atop uv\in E(G)}
\phi(uv)-\sum_{u\in V(G)\atop vu\in E(G)} \phi(vu).
\]
So $\phi$ is an $S$-flow if an only if $\partial_\phi=0$. Consider
the expression
\[
A=\sum_{\phi:~E(G)\to S} \prod_{v\in V(G)}
\sum_{\chi\in\Gamma^*}\chi(\partial_\phi(v)).
\]
The summation over $\chi$ is $0$ unless $\partial_\phi(v)=0$, in
which case it is $m$. So the product over $v\in V(G)$ is $0$ unless
$\phi$ is an $S$-flow, in which case it is $m^n$. So $A\cdot m^{-n}$
counts $S$-flows.

On the other hand, we can expand the product over $v\in V(G)$; each
term will correspond to a choice of a character $\psi_v$ for each
$v$, and so we get
\[
A=\sum_{\phi:~E(G)\to S} \sum_{\psi:~V(G)\to\Gamma^*}\prod_{v\in
V(G)} \psi_v(\partial_\phi(v)).
\]
Here (using that $\psi_v$ is a character)
\[
\psi_v(\partial_\phi(v))=\prod_{u\in V(G)\atop uv\in E(G)}
\psi_v(\phi(uv))\prod_{u\in V(G)\atop vu\in E(G)}
\psi_v(\phi(vu))^{-1},
\]
so we get that
\[
A=\sum_{\phi:~E(G)\to S} \sum_{\psi:~V(G)\to\Gamma^*} \prod_{uv\in
E(G)} \psi_v(\phi(uv)) \psi_u(\phi(uv))^{-1}.
\]
Interchanging the summation, the inner sum factors:
\begin{align*}
\sum_{\phi:~E(G)\to S} \prod_{uv\in E(G)} \psi_v(\phi(uv))
\psi_u(\phi(uv))^{-1}&=\prod_{uv\in E(G)} \sum_{s\in S} \psi_v(s)
\psi_u(s)^{-1}\\
&=\prod_{uv\in E(G)} \beta_{\psi_u\psi_v}=m^n w(\psi),
\end{align*}
showing that
\[
A=m^n \hom(G,H).
\]
This proves (\ref{FLOWHOM}).

\smallskip

\end{example}

\begin{example}[The role of multiple edges]\label{MULTIPLE}
Let us give an example of a reflection positive graph parameter $f$
for which $r_f(k)$ is finite for every $k$ but has superexponential
growth. The example also shows that we have to be careful with
multiple edges. Let, for each graph $G$, $G'$ denote the graph which
we obtain from $G$ by keeping only one copy of each parallel class of
edges. Let
\[
f(G)=2^{-|E(G')|}.
\]
It is not hard to see that the connection matrices $M(f,k)$ are
positive semidefinite. This graph parameter is, in fact, the limit of
parameters of the form $\hom(.,H)$: take homomomorphisms into a
random graph $H=G(n,1/2)$, with all node weights $=1/n$ and all
edge-weights 1. Furthermore, it is not hard to see that the rank of
$M(f,k)$ is $2^{k\choose 2}$. This is finite but superexponential, so
the parameter is not of the form $\hom(.,H)$.

Note, however, that for a simple graph $G$ ({\it i.e.}, if $G$ has no
multiple edges), $f(G)=2^{|E(G)|}$ can be represented as the number
of homomorphisms into the graph consisting of a single node with a
loop, where the node has weight 1 and the loop has weight $1/2$.
\end{example}

\begin{example}[Chromatic polynomial]\label{TUTTE}
The following example is from \cite{FLW}. Let $p(G)=p(G,x)$ denote
the chromatic polynomial of the graph $G$. For every fixed $x$, this
is a multiplicative graph parameter. To describe its
rank-connectivity, we need the following notation. For $k,q\in\Z_+$,
let $B_{kq}$ denote the number of partitions of a $k$-element set
into at most $q$ parts. So $B_k=B_{kk}$ is the $k$-th Bell number.
With this notation,
\[
r_p(k)=
  \begin{cases}
    B_{kx} & \text{if $x$ is a nonnegative integer}, \\
    B_k & \text{otherwise}.
  \end{cases}
\]
Note that this is always finite, but if $x\notin \Z_+$, then it grows
faster than $c^k$ for every $c$.

Furthermore, $M(p,k)$ is positive semidefinite if and only if either
$x$ is a positive integer or $k\le x+1$. The parameter $M(p,k)$ is
reflection positive if and only if this holds for every $k$, i.e., if
and only if $x$ is a nonnegative integer, in which case indeed
$p(G,x)=\hom(G,K_x)$.
\end{example}

\begin{example}[Homomorphisms into infinite graphs]\label{INFINITE}
We can extend the definition of $\hom(G,H)$ to infinite weighted
graphs $H$ provided the node and edge-weights form sufficiently fast
convergent sequences. Reflection positivity remains valid, but the
rank of $M(f_H,k)$ will become infinite, so this graph parameter
cannot be represented by a finite $H$.

More generally, let $a>0$, $I=[0,a]$, and let $W:~I\times I\to\R$ be
a measurable function such that for every $n\in\Z_+$,
\[
\int_0^a \int_0^a |W(x,y)|^n\,dx\,dy<\infty.
\]
Then we can define a graph parameter $f_W$ as follows. Let $G$ be a
finite graph on $n$ nodes, then
\[
f_W(G)=\int_{I^n} \prod_{ij\in E(G)} W(x_i,x_j)\,dx_1\,\dots\,dx_n.
\]

It is easy to see that for every weighted (finite or infinite) graph
$H$, the graph parameter $f_H$ is a special case. Furthermore, $f_W$
is reflection positive. However, it can be shown that the graph
parameter in Example \ref{MULTIPLE} cannot be represented in this
form \cite{SZ}.
\end{example}

\section{Proof of Theorem \ref{MAIN}}

\subsection{The algebra of graphs}

In this first part of the proof, we only assume that for every $k$,
$M(f,k)$ is positive semidefinite, has finite rank $r_k$, and
$r_0=1$. We know that $f$ is multiplicative, and hence $f(K_0)=1$. We
can replace the parameter $f(G)$ by $f(G)/f(K_1)^{-|V(G)|}$; this can
be reversed by scaling of the node-weights of the target graph $H$
representing $f$, once we have it constructed. So we may assume that
$f(K_1)=1$. Combined with multiplicativity, this implies that we can
delete (or add) isolated nodes from any graph $G$ without changing
$f(G)$.

It will be convenient to put all $k$-labeled graphs into a single
structure as follows. By a {\it partially labeled graph} we mean a
finite graph in which some of the  nodes are labeled by distinct
nonnegative integers. Two partially labeled graphs are {\it
isomorphic} if there is an isomorphism between them preserving all
labels. For two partially labeled graphs $G_1$ and $G_2$, let
$G_1G_2$ denote the partially labeled graph obtained by taking the
disjoint union of $G_1$ and $G_2$, and identifying nodes with the
same label. For every finite set $S\subseteq\Z_+$, we call a
partially labeled graph {\it $S$-labeled}, if its labels form the set
$S$.

Let $\GG$ denote the (infinite dimensional) vector space of formal
linear combinations (with real coefficients) of partially labeled
graphs. We can turn $\GG$ into an algebra by using $G_1G_2$
introduced above as the product of two generators, and then extending
this multiplication to the other elements linearly. Clearly $\GG$ is
associative and commutative, and the empty graph is a unit element.

For every finite set $S\subseteq \Z_+$, the set of all formal linear
combinations of $S$-labeled graphs forms a subalgebra $\GG(S)$ of
$\GG$. The graph with $|S|$ nodes labeled by $S$ and no edges is a
unit in this algebra, which we denote by $U_S$.

We can extend $f$ to a linear functional on $\GG$, and define an
inner product
\[
\langle x,y\rangle = f(xy)
\]
for $x,y\in\GG$. By our hypothesis that $f$ is reflection positive,
this inner product is positive semidefinite, i.e., $\langle
x,x\rangle\ge 0$ for all $x\in\GG$. Indeed, if $x=\sum_{G\in\GG}x_G
G\in\GG$ (with a finite number of nonzero terms), then
\[
\langle x,x\rangle =
f(xx)=\sum_{G_1,G_2\in\GG}x_{G_1}x_{G_2}f(G_1G_2) \ge 0,
\]
since the quadratic form is of the form $x\T M(f,k)x$ for a large
enough $k$, and $M(f,k)$ is positive semidefinite.

Let
\[
\KK= \{x\in\GG: ~f(xy)=0 ~\forall y\in\GG\}
\]
be the annihilator of $\GG$. It follows from the positive
semidefiniteness of the inner product that $x\in \KK$ could also be
characterized by $f(xx)=0$.

Clearly, $\KK$ is an ideal in $\GG$, so we can form the quotient
algebra $\hat\GG=\GG/\KK$. We can also define
$\hat\GG(S)=\GG(S)/\KK$. It is easy to check that every graph $U_S$
has the same image under this factorization, namely unit element $u$
of $\hat\GG$. Furthermore, if $x\in \KK$ then $f(x)=f(xu)=0$, and so
$f$ can also be considered as a linear functional on $\hat\GG$. We
denote by $\hat G$ the element of $\hat\GG$ corresponding to the
partially labeled graph $G$.

Next, note that $\hat\GG(S)$ is a finite dimensional commutative
algebra of dimension $r_{|S|}$, with the positive definite inner
product $\langle x,y\rangle = f(xy)$.

Let $S\subseteq T$ be finite subsets of $\Z_+$. Then
$\hat\GG(S)\subseteq\hat\GG(T)$. Indeed, every $S$-labeled graph $G$
can be turned into a $T$-labeled graph $G'$ by adding $|T\setminus
S|$ new isolated nodes, and label them by the elements of $T\setminus
S$. It is straightforward to check that $G-G'\in\KK$, and so $G$ and
$G'$ correspond to the same element of $\hat\GG$.

We'll also need the orthogonal projection $\pi_S$ of $\hat\GG$ to the
subalgebra $\hat\GG(S)$. This has a very simple combinatorial
description. For every partially labeled graph $G$ and
$S\subseteq\Z_+$, let $G_S$ denote the $S$-labeled graph obtained by
deleting the labels not in $S$; then $\pi_S(\hat{G})=\hat{G}_S$.

To see this, let $G$ be any partially labeled graph. Then
\[
\langle G_S, G-G_S\rangle = f(G_S(G-G_S))=f(G_SG)-f(G_SG_S)=0,
\]
since $GG_S$ and $G_SG_S$ are isomorphic graphs. Hence
\[
\langle \hat{G}_S, \hat{G}-\hat{G}_S\rangle =0,
\]
showing that $\hat{G}_S$ is indeed the orthogonal projection of
$\hat{G}$ onto $\hat\GG(S)$.

We'll be interested in the idempotent elements of $\hat\GG$. If $p$
is idempotent, then
\[
f(p)=f(pp)=\langle p,p\rangle >0.
\]
For two idempotents $p$ and $q$, we say that $q$ {\it resolves} $p$,
if $pq=q$. It is clear that this relation is transitive.

Let $S$ be a finite subset of $\Z_+$, and set $r=r_{|S|}$. Since the
algebra $\hat\GG(S)$ is finite dimensional and commutative, and all
its elements are self-adjoint with respect to the positive definite
inner product $\langle.,.\rangle$, it has a (uniquely determined)
basis $\PP_S=\{p^S_1,\dots, p^S_r\}$ such that $(p^S_i)^2=p^S_i$ and
$p^S_ip^S_j=0$ for $i\not=j$. We denote by $\PP_{T,p}$ the set of all
idempotents in $\PP_T$ that resolve a given idempotent $p$. If
$|T|=|S|+1$, then the number of elements in $\PP_{T,p}$ will be
called the {\it degree} of $p\in\PP_S$, and denoted by $\deg(p)$.
Obviously this value is independent of which ($|S|+1$)-element
superset $T$ of $S$ we are considering.

\begin{claim}\label{RESOLVESUM}
Let $r$ be any idempotent element of $\GG(S)$. Then $r$ is the
sum of those idempotents in $\PP_S$ that resolve it.
\end{claim}

Indeed, we can write
\[
r=\sum_{p\in\PP_S} \mu_p p
\]
with some scalars $\mu_p$. Now using that $r$ is idempotent:
\[
r=r^2=\sum_{p,p'\in\PP_S} \mu_p\mu_{p'} pp'= \sum_{p\in\PP_S} \mu_p^2
p,
\]
which shows that $\mu_p^2=\mu_p$ for every $p$, and so
$\mu_p\in\{0,1\}$. So $r$ is the sum of some subset
$X\subseteq\PP_S$. It is clear that $rp=p$ for $p\in X$ and $rp=0$
for $p\in \PP_S\setminus X$, so $X$ consists of exactly those
elements of $\PP_S$ that resolve $q$.

As a special case, we see that
\begin{equation}\label{UNIT}
u=\sum_{p\in \PP_S} p
\end{equation}
is the unit element of $\hat\GG(S)$ (this is the image of the graph
$U_S$), and also the unit element of the whole algebra $\hat\GG$.

\begin{claim}\label{RESOLVEPART}
Let $S\subset T$ be two finite sets. Then every $q\in\PP_T$ resolves
exactly one element of $\PP_S$.
\end{claim}

Indeed, we have by (\ref{UNIT}) that
\[
u = \sum_{p\in\PP_S}p=\sum_{p\in\PP_S}\sum_{q\in\PP_T\atop q\text{
resolves }p} q,
\]
and also
\[
u =\sum_{q\in\PP_T} q,
\]
so by the uniqueness of the representation we get that every $q$ must
resolve exactly one $p$.

\begin{claim}\label{UNLABEL}

Let $S,T,U$ be finite sets and $S=T\cap U$. If $x\in\hat\GG(T)$ and
$y\in\hat\GG(U)$, then
\[
f(xy)=f(\pi_S(x)y).
\]
\end{claim}

Indeed, for every $T$-labeled graph $G_1$ and $U$-labeled graph
$G_2$, the graphs $G_1G_2$ and $\pi_S(G_1)G_2$ are isomorphic. Hence
the claim follows by linearity.

\begin{claim}\label{RESOLVEPROJ}
If $p\in\PP_S$ and $q$ resolves $p$, then
\[
\pi_S(q)=\frac{f(q)}{f(p)} p.
\]
\end{claim}

We show that both sides give the same inner product with every basis
element in $\PP_S$. Since $q$ does not resolve any $p'\in
\PP_S\setminus \{p\}$, we have $p'q=0$ for every such $p'$. By Claim
\ref{UNLABEL}, this implies that
\[
\langle p',\pi_S(q)\rangle = f(p'\pi_S(q))=f(p'q)=0 =\langle
p',\frac{f(q)}{f(p)} p\rangle.
\]
Furthermore,
\[
\langle p,\pi_S(q)\rangle = f(p\pi_S(q))=f(pq)=f(q) =\langle
p,\frac{f(q)}{f(p)} p\rangle.
\]
This proves the Claim.

\begin{claim}\label{RESOLVEF}
Let $S,T,U$ be finite sets and $S=T\cap U$. Then for any $p\in
\PP_S$, $q\in \PP_{T,p}$ and $r\in \hat\GG(U)$ we have
\[
f(p)f(qr)=f(q)f(pr).
\]
\end{claim}

Indeed, by Claims \ref{RESOLVEPROJ} and \ref{UNLABEL},
\[
f(qr) = f(\pi_S(q)r) = \frac{f(q)}{f(p)}f(pr).
\]

\begin{claim}\label{RESOLVEPROD}
If both $q\in \PP_T$ and $r\in \PP_U$ resolve $p$, then $qr\not=0$.
\end{claim}

Indeed, by Claim \ref{RESOLVEF},
\[
f(qr)=\frac{f(q)}{f(p)}f(pr)=\frac{f(q)}{f(p)}f(r) >0.
\]

\begin{claim}\label{RESOLVEDE(G)}
If $S\subset T$, and $q\in\PP_T$ resolves $p\in\PP_S$, then
$\deg(q)\ge\deg(p)$.
\end{claim}

It suffices to show this in the case when $|T|=|S|+1$. Let
$U\subset\Z_+$ be any $(|S|+1)$-element superset of $S$ different
from $T$. Let $Y$ be the set of elements in $\PP_U$ resolving $p$.
Then $p=\sum_{r\in Y} r$ by Claim \ref{RESOLVESUM}. Here
$|Y|=\deg(p)$. Furthermore, we have
\[
\sum_{r\in Y} rq = q\sum_{r\in Y} r = qp=q.
\]
Each of the terms on the left hand side is nonzero by Claim
\ref{RESOLVEPROD}, and since the terms are all idempotent, each of
them is the sum of one or more elements of $\PP_{T\cup U}$.
Furthermore, if $r,r'\in Y$ ($r\not= r'$), then we have the
orthogonality relation
\[
(rq)(r'q)=(rr')q=0,
\]
so the basic idempotents in the expansion of each term are different.
So the expansion of $q$ in $\PP_{T\cup U}$ contains at least
$|Y|=\deg(p)$ terms, which we wanted to prove.

\subsection{Bounding the expansion}

From now on, we assume that there is a $q>0$ such that $r_k\le q^k$
for all $k$.

So if a basic idempotent $p\in\PP_S$ has degree $D$, then there are
$D$ basic idempotents on the next level with degree $\ge D$, and
hence if $|T|\ge |S|$, then the dimension of $\hat\GG(T)$ is at least
$D^{|T\setminus S|}$. It follows that the degrees of basic
idempotents are bounded by $q$; let $D$ denote the maximum degree,
attained by some $p\in\PP_S$.

Let us fix such a set $S$ and $p\in\PP_S$ with maximum degree $D$.
For $u\in\Z_+\setminus S$, let $q_1^u,\dots, q_D^u$ denote the
elements of $\PP_{S\cup\{u\}}$ resolving $p$. Note that for
$u,v\in\Z_+\setminus S$, there is a natural isomorphism between
$\hat\GG(S\cup\{u\})$ and $\hat\GG(S\cup\{v\})$ (induced by the map
that fixes $S$ and maps $u$ onto $v$), and we may choose the labeling
so that $q_i^u$ corresponds to $q_i^v$ under this isomorphism.

Next we describe, for a finite set $T\supset S$, all basic
idempotents in $\PP_T$ that resolve $p$. Let $V=T\setminus S$, and
for every map $\phi:~V\to \{1,\dots,D\}$, let
\begin{equation}\label{QPHI-DEF}
q_\phi=\prod_{v\in V} q^v_{\phi(v)}.
\end{equation}
Note that by Claim \ref{RESOLVEF},
\begin{equation}\label{21}
f(q_{\phi})= f(\prod_{v\in V}q^v_{\phi(v)})= \bigl(\prod_{v\in
V}\frac{f(q^v_{\phi(v)})}{f(p)}\bigr)f(p) \neq 0,
\end{equation}
and so $q_{\phi}\neq 0$.

\begin{claim}\label{RESOLVEALL}
\[
\PP_{T,p}=\{q_{\phi}:~\phi\in \{1,\ldots,D\}^V\}.
\]
\end{claim}

We prove this by induction on $|T\setminus S|$. For $|T\setminus
S|=1$ the assertion is trivial. Suppose that $|T\setminus S|>1$. Let
$u\in T\setminus S$, $U=S\cup \{u\}$ and $W=T\setminus \{u\}$; thus
$U\cap W=S$. By the induction hypothesis, the basic idempotents in
$\PP_W$ resolving $p$ are elements of the form $q_\psi$
($\psi\in\{1,\dots,D\}^{V\setminus \{u\}}$).

Let $r$ be one of these. By Claim \ref{RESOLVEPROD}, $rq_i^u\not=0$
for any $1\le i\le D$, and clearly resolves $r$. We can write
$rq_i^u$ as the sum of basic idempotents in $\PP_T$, and it is easy
to see that these also resolve $r$. Furthermore, the basic
idempotents occurring in the expression of $rq^u_i$ and $rq^u_j$
($i\not= j$) are different. But $r$ has degree $D$, so each $rq_i^u$
must be a basic idempotent in $\PP_T$ itself.

Since the sum of the basic idempotents $r q_i^u$  ($r\in \PP_{W,p}$,
$1\le i\le D$) is $p$, it follows that these are all the elements of
$\PP_{T,p}$. This proves the Claim.

It is immediate from the definition that an idempotent $q_\phi$
resolves $q_i^v$ if and only $\phi(v)=i$. Hence it also follows that
\begin{equation}\label{QISUM}
q^v_{i}=\sum_{\phi:~\phi(v)=i}q_{\phi}.
\end{equation}

\subsection{Constructing the target graph}

Now we can define $H$ as follows. Let $H$ be the looped complete
graph on $V(H)=\{1,\dots,D\}$. We have to define the node weights and
edge weights.

Fix any $u\in \Z_+\setminus S$. For every $i\in V(H)$, let
\[
\alpha_i=\frac{f(q_i^u)}{f(p)}
\]
be the weight of the node $j$. Clearly $\alpha_i>0$.

Let $u,v\in \Z_+\setminus S$, $v\not=u$, and let $W=S\cup\{u,v\}$.
Let $K_{uv}$ denote the graph on $W$ which has only one edge
connecting $u$ and $v$, and let $k_{uv}$ denote the corresponding
element of $\hat\GG(W)$. We can express $pk_{uv}$ as a linear
combination of elements of $\PP_{W,p}$ (since for any
$r\in\PP_W\setminus \PP_{W,p}$ one has $rp=0$ and hence
$rpk_{u,v}=0$):
\[
pk_{uv}=\sum_{i,j}\beta_{ij}q^u_iq^v_j.
\]
This defines the weight $\beta_{ij}$ of the edge $ij$. Note
$\beta_{ij}=\beta_{ji}$ for all $i,j$, since $pk_{uv}=pk_{vu}$.

We prove that this weighted graph $H$ graph gives the right
homomorphism function.

\begin{claim}\label{FHOM}
For every finite graph $G$, $f(G)=\hom(G,H)$.
\end{claim}

By (\ref{QISUM}), we have for each pair $u,v$ of distinct elements of
$V(G)$
\[
pk_{uv}= \sum_{i,j\in V(H)}\beta_{i,j}q^u_{i}q^v_{j}= \sum_{i,j\in
V(H)}\beta_{i,j}\sum_{\phi:~\phi(u)=i\atop ~~~\phi(v)=j}q_{\phi}
=\sum_{\phi\in V(H)^V}\beta_{\phi(u),\phi(v)}q_{\phi}.
\]

Consider now any $V$-labeled graph $G$ with $V(G)=V$, and let $g$ be
the corresponding element of $\hat\GG$. Then
\begin{align*}
pg&=\prod_{uv\in E(G)}pk_{uv} =\prod_{uv\in E(G)}\bigl(\sum_{\phi\in
V(H)^V}\beta_{\phi(u),\phi(v)}q_{\phi}\bigr)\\
& =\sum_{\phi:V\to V(H)}\bigl(\prod_{uv\in
E(G)}\beta_{\phi(u),\phi(v)}\bigr) q_{\phi}.
\end{align*}
Since $p\in \hat\GG(S)$, $g\in \hat\GG(V)$ and $S\cap V=\emptyset$,
we have $f(p)f(g)=f(pg)$, and so by (\ref{21}),
\begin{align*}
f(p)f(g)&=f(pg) =\sum_{\phi\in V(H)^V}\bigl(\prod_{uv\in E(G)}
\beta_{\phi(u),\phi(v)}\bigr)f(q_{\phi})\\
&=\sum_{\phi:V\to V(H)}\bigl(\prod_{uv\in E(G)}
\beta_{\phi(u),\phi(v)}\bigr)\bigl(\prod_{v\in V(G)}
\alpha_{\phi(v)}\bigr)f(p),
\end{align*}
The factor $f(p)>0$ can be cancelled from both sides, completing the
proof.

\section{Extensions: Graphs with loops, directed graphs and hypergraphs}

In our arguments we allowed parallel edges in $G$, but no loops.
Indeed, the representation theorem is false if $G$ can have loops: it
is not hard to check that the graph parameter
\[
f(G)=2^{\#\text{loops}}
\]
cannot be represented as a homomorphism function, even though its
connection matrix $M(f,k)$ is positive semidefinite and has rank 1.
To get a representation theorem for graphs with loops, each loop $e$
in the target graph $H$ must have two weights: one which is used when
a non-loop edge of $G$ is mapped onto $e$, and the other, when a loop
of $G$ is mapped onto $e$. With this modification, the proof goes
through.

The constructions and results above are in fact more general; they
extend to directed graphs and hypergraphs. There is a common
generalization to these results, using semigroups. This will be
stated and proved in a separate paper.

\section*{Acknowledgement}

We are indebted to Christian Borgs, Jennifer Chayes, Monique Laurent,
Miki Simonovits, Vera T.~S\'os, Bal\'azs Szegedy, G\'abor Tardos and
Kati Vesztergombi for many valuable discussions and suggestions on
the topic of graph homomorphisms.

\end{document}